\documentclass[a4paper,11pt]{article}
\usepackage[english]{babel}
\usepackage[psamsfonts]{amssymb}
\usepackage{cmmib57}
\usepackage{exscale}
\usepackage{amsmath}
\usepackage{amscd}
\usepackage{latexsym}
\numberwithin{equation}{section}
\newtheorem{theorem}{Theorem}
\newtheorem{corollary}{Corollary}
\newtheorem{lemma}{Lemma}

\newtheorem{proposition}{Proposition}

\newcommand{\cqd}{\hfill$\Box$}
\newcommand{\R}{{\mathbb R}}

\newcommand{\ind}{\mbox{\Large$\chi$}}
\newcommand{\tor}{{\mathbb T}}
\newcommand{\G}{{\mathbb G}}

\begin{document}
\begin{titlepage}
\mbox{}
\begin{center}
{\bf \Large A lattice scheme for\\
stochastic partial differential equations\\
of elliptic type in dimension $d\ge 4$}

\vspace*{2cm}
\small
\begin{tabular}{l@{\hspace{2cm}}l}
{\small\sc Teresa Mart\'{\i}nez}$^{\dagger}$&{\small\sc Marta Sanz-Sol\'e}$^{\,*}$\\
Departamento de Matem\'aticas & Facultat de Matem\`atiques\\
Universidad Aut\'onoma de Madrid,&Universitat de Barcelona\\
Campus de Cantoblanco, &Gran Via 585\\
E-28049 Madrid &E-08007 Barcelona\\ 
teresa.martinez@uam.es& marta.sanz@ub.edu\\
\end{tabular}
\end{center}

\vspace*{3cm}

\noindent{\bf Abstract:} We study a stochastic boundary value problem on $(0,1)^d$ of elliptic type in dimension $d\ge 4$, driven by a coloured noise.  An approximation scheme based on a suitable discretization of the Laplacian on a lattice of 
$(0,1)^d$ is presented; we also give the rate of convergence to the original SPDE in $L^p(\Omega;L^{2}(D))$--norm, for some values of $p$. 
\bigskip

\smallskip


\vfill

\begin{itemize}
\item[{$^\dagger$}]{\small Partially supported by the grant BFM 2002-04013-C02-02 from the
Direcci\'on General de
Universidades, Ministerio de Educaci\'on y Ciencia, Spain.}
\item[{$^*$}]{\small Supported by the grant BFM 2003-01345 from the Direcci\'on
General de Investigaci\'on, Ministerio de Educaci\'on y Ciencia, Spain.}\\[-.8cm]
\end{itemize}
\end{titlepage}
\date{\today}

\section{Introduction}\label{s1}
We consider the stochastic boundary value problem of elliptic type
\begin{eqnarray}\label{1.1}
\Delta u(x)-f(u(x))&=&g(x)+\dot F(x),\quad\mbox{ for }x\in D\\
u(x)&=&0,\quad\mbox{ for }x\in \partial D\nonumber,
\end{eqnarray}
where $D=(0,1)^d$, $d\geq 4$, and $\dot F(x)$
is the coloured noise defined at the beginning of Section \ref{s2}. The coefficient $g$ belongs to $L^2(D)$ and $f$ is of the form
\begin{equation}\label{1.2}
f(x)=f_1(x)+f_2(x), 
\end{equation}
with  $f,f_1,f_2:\R\rightarrow\R$ satisfying\\
\noindent$(f1)$\quad $f_1$ is continuous, non-decreasing and bounded, $\sup_{x\in \R}|f_1(x)|\le M$, \\
\noindent$(f2)$\quad $f_2$ is Lipschitz with {\it small} Lipschitz
constant $L$.\\

Let  $G_D$ be  the Green function of the Poisson equation
$\Delta v = b$, with boundary condition $v(x)= 0$  if $x\in\partial D$. By a solution of (\ref{1.1}) we mean a stochastic process
$(u(x), x\in D)$ satisfying
\begin{equation}
\label{1.3}
u(x)=\int_DG_D(x,y)f(u(y))\,dy+\int_DG_D(x,y)g(y)\,dy+\int_DG_D(x,y)dF(y).
\end{equation}

For dimensions $d=1,2,3$, $\dot F$ a white noise and $D$ a bounded domain in $\R^d$ with regular boundary, 
the existence of a unique solution to equation (\ref{1.3}) has been proved in  \cite{bu-pa:90} (see also \cite{catherine}).
In this framework, for the particular domain $D=(0,1)^d$, Gy\"ongy and Mart\'{\i}nez introduced in \cite{GM} a numerical scheme based on discretization of the Laplacian
and  gave the rate of convergence in the $L^2(D)$--norm (see also \cite{al-no-zh:98} for related work).
For these range of dimensions the Green function satisfies $\sup_{x\in D}\Vert G_d(x,\cdot)\Vert_{L^2(D)}<\infty$. This fact has two important consequences.
Firstly,  $\left(\int_DG_D(x,y)dF(y), x\in D\right)$ defines a Gaussian process; secondly, one can work with the Fourier series of $G_D$. Actually,
the lattice scheme in \cite{GM} is related to truncated Fourier expansions.

For parabolic equations driven by space-time white noise, numerical schemes based on lattice approximations have been introduced in \cite{gy:98} (see also \cite{gy:99}). The method set up in these papers has been successfully applied in \cite{GM} and also to other types of SPDEs, for instance, to a wave equation in spatial dimension 1 in \cite{qss}.
 
In this paper, we study a lattice approximation for the stochastic partial differential equation (\ref{1.1}) for $d\ge 4$, extending the results  of \cite{GM}. 

For dimensions $d\ge 3$, the Green function associated with (\ref{1.1}) is given by
\begin{equation}
\label{1.4}
G_D(x,y)=G(x,y)+E_x(G(B_\tau,y)),\quad\mbox{with}\quad G(x,y)=\frac{C_d}{|x-y|^{d-2}}.
\end{equation}
In these expresions, $C_d$ is a constant and $B$ is a $d$-dimensional Brownian motion starting at $x$, stopped at $\tau$
-its first exit time of $D$ (see \cite{doob} for details). 

We will prove that for $d\ge 4$, $G_D(x,\cdot)\in L^\alpha(D)$, $\alpha\in[1,d/(d-2))$, uniformly in
$x$.  Hence, we cannot use a $L^2$ theory. In particular, the stochastic integral $\int_DG_D(x,y)dF(y)$ with respect to a white noise 
cannot be defined as a real-valued $L^2$ random variable. This problem forces the choice of a coloured noise, as a way to give a rigurous meaning to \eqref{1.1}.

The contents of the paper are as follows. Section \ref{s2} is devoted to study the spde. First, we attach a precise meaning to the
stochastic integral term in (\ref{1.3}) as a Gaussian process indexed by elements of its reproducing kernel Hilbert space and
give sample path properties related with the regularity of the covariance of the noise. Secondly, we prove existence and uniqueness
of solution. The methods are common to those of nonlinear monotone operator equations (see \cite{lions}, \cite{ze:90}) and those
used also in \cite{{bu-pa:90}} and \cite{catherine}. However, since our setting is different, we feel interesting to give the details.
We also analyze properties of the solution; in particular, H\"older continuity of the sample paths. Section \ref{s3} is devoted to the
lattice approximation and the analysis of the rate of convergence (see Theorem \ref{t3.19}). First, we introduce a family of smoothed Green functions
obtained by convolution of $G_D(x,\cdot)$ with an approximation of the identity. Then, we introduce a lattice approximation which corresponds
to a weighted discretization of the Laplacian on a given grid, with weights related to the smoothing. The approximation result follows from a careful
analysis of the errors in both of these approximations, taking as smoothing parameter an appropriate function of the norm of the grid.

Let us give some indications about notation.  By $|x|$ we denote the Hilbert-Schmidt norm of any element
$x\in\R^d$. The letter $C$ denotes, unless otherwise stated, a constant that may not be the
same from one occurrence to another. Sometimes, we denote fixed values of
constants by adding a subindex, and the dependence on some parameters (as the dimension) with an
argument, e.g. $C_1$ or $C(d)$. Throughout the work, the symbols $\beta=(\beta_1,\dots,\beta_d)$,
$i=(i_1,\dots,i_d)$ denote indexes belonging to the sets
$I^d=\{1,2,\dots\}^d$, $I^d_n=\{1,\dots,n-1\}^d$,
with $|\beta|^2=\beta_1^2+\dots+\beta_d^2$. Observe that for
$\beta\in I^d\backslash I^d_n$, $|\beta|^2\geq n^2$, and for $\beta\in I^d_n$, $d\leq
|\beta|^2\leq dn^2$. 

\section{Study of the equation}
\label{s2}
Let $\varphi$ be the density of a non negative measure on $\R^d$, non negative definite and tempered.
We consider a centered Gaussian process $F(\psi)$, indexed by the space $\mathcal{D}(\R^d)$ of Schwartz test functions,
with covariance function
\begin{equation}\label{2.1}
E(F(\psi_1),F(\psi_2))=\int_{\R^d}\int_{\R^d}\psi_1(x)\varphi(x-y)\psi_2(y)\,dx\,dy,
\end{equation}
$\psi_1,\psi_2\in\mathcal{D}(\R^d)$, defined on some probability space 
$(\Omega,{\mathcal F},P)$.

Let $\mathcal{H}$ denote the completion of the inner product space consisting of functions $\psi\in\mathcal{D}(\R^d)$
endowed with the inner product
\begin{equation}
\label{2.2}
\langle \psi_1,\psi_2\rangle_{\mathcal{H}}=\int_{\R^d}\int_{\R^d}\psi_1(x)\varphi(x-y)\psi_2(y)\,dx\,dy.
\end{equation}
The space $\mathcal{H}$ is the reproducing kernel Hilbert space corresponding to $F$.

We consider the stochastic
partial differential equation  \eqref{1.1}, with the assumptions given in Section \ref{s1}.
We notice for further use the following property: 

\noindent{\bf(P)}\, $f$ is a function of the form (\ref{1.2}) with $f_1$ non-decreasing
and $f_2$ Lipschitz with Lipschitz constant $L$, if and only if
for every $u,v\in\R$, 
\begin{equation}\label{2.3} 
(u-v)(f(u)-f(v))\geq -L(u-v)^2. 
\end{equation}


\subsection{The stochastic integral}\label{s2.1}

In order to give a rigourous meaning to \eqref{1.1}, we have to precise what type of stochastic integral we are considering. By classical results on abstract Wiener spaces, $h\in {\mathcal H}\mapsto F(h)$ defines a linear continuous
functional that satisfies \eqref{2.1} (see  \cite{nualart}). Thus, the stochastic integral term in \eqref{1.3} is well defined as long as we prove that $G_D(x,\cdot)\in{\mathcal H}$. The next lemma provides a sufficient condition for this property to hold.

\begin{lemma}\label{l2.5} Let $p\in[1,\infty)$, $\frac{1}{p'}+\frac{1}{p}=1$.
\begin{enumerate}
\item{(1)} Let  $\psi\in L^p(\R^d)$, $\varphi\in L^{\frac{p'}{2}}(\R^d)$. Then,
\begin{equation*}
\left\vert\int_{\R^d}\int_{\R^d}\psi(x)\varphi(x-y)\psi(y)\,dx\,dy\right\vert
\leq \|\psi\|_{L^p(\R^d)}^2\|\varphi\|_{L^{\frac{p'}{2}}(\R^d)}.
\end{equation*}
\item{(2)} Assume $\psi\in L^1(\R^d)\cap L^p(\R^d)$, $\varphi\in L^{p'}(\R^d)$. Then,
\begin{equation*}
\left\vert\int_{\R^d}\int_{\R^d}\psi(x)\varphi(x-y)\psi(y)\,dx\,dy\right\vert
\le \|\psi\|_{L^p(\R^d)}\|\psi\|_{L^1(\R^d)}\|\varphi\|_{L^{p'}(\R^d)}.
\end{equation*}
\end{enumerate}
\end{lemma}

\noindent{\sc Proof}. Applying first H\"older's inequality and then Young's inequality for convolutions (\cite{adams}, Corollary 2.25)
yield
\begin{eqnarray*}
&&\bigg|\int_{\R^d}\int_{\R^d}\psi(x)\varphi(x-y)\psi(y)\,dx\,dy\bigg|=\bigg|\int_{\R^d}\psi(x)(\varphi*\psi)(x)\,dx\bigg|\\
&&\leq\|\psi\|_{L^p(\R^d)}\|\varphi*\psi\|_{L^{p'}(\R^d)}\leq \|\psi\|_{L^p(\R^d)}^2\|\varphi\|_{L^{p'/2}(\R^d)}.
\end{eqnarray*}
Young's theorem (\cite{adams}, Theorem 2.24) implies
\begin{equation*}
\bigg|\int_{\R^d}\psi(x)(\varphi*\psi)(x)\,dx\bigg|
\le \|\psi\|_{L^1(\R^d)}\|\varphi\|_{L^{p'}(\R^d)}\|\psi\|_{L^p(\R^d)}.
\end{equation*}
Both estimates yield the lemma.
\cqd

We next prove a basic result on the Green function $G_D$.
\begin{lemma}\label{l2.6}
For any $\alpha\in [1,\frac{d}{d-2})$, there exists a positive constant $C_1$ depending on $\alpha$ and $d$, such that
$$
\sup_{x\in D}\|G_D(x,\cdot)\|_{L^\alpha(D)}\leq C_1.
$$
Consequently, $\|G_D\|_{L^\alpha(D\times D)}\leq C$.
\end{lemma}

\noindent{\sc Proof}. Set  $I_1=\Vert G(x,\cdot)\Vert_{L^\alpha(D)}$. Clearly, 
$$
I^\alpha_1=C_d^\alpha \int_{x-D}\frac{dz}{|z|^{(d-2)\alpha}}
\leq C(d)\int_0^2\, r^{d(1-\alpha)-1+2\alpha}\,dr.
$$
The last integral is finite if and only if $\alpha<\frac{d}{d-2}$; in this case, its value is a constant $C(\alpha,d)$
independent of $x$.

Let $I_2=\|E_x(G(B_\tau,y))\|_{L^\alpha(D)}$ and let $P^x$  be the law of the random variable $B_{\tau}$.
By Minkowski inequality,
\begin{equation*}
I_2=\bigg\|\int_{\R^d} G(u,\cdot)\,dP^x(u)\bigg\|_{L^\alpha(D)}
\leq \int_{\R^d}\|G(u,\cdot)\|_{L^\alpha(D)}\,dP^x(u)\leq \tilde C(\alpha,d).
\end{equation*}

Since by \eqref{1.4}, $\|G_D(x,\cdot)\|_{L^\alpha(D)}\leq I_1+I_2$, with the upper bounds obtained so far, we finish the proof of the lemma.

 \cqd

Fix $\alpha\in[1,\frac{d}{d-2})$ and denote by $\alpha'$ its conjugate. 
Set ${\mathcal L}^\alpha=L^{\alpha'}(\R^d)\cup L^{\alpha'/2}(\R^d)$.
 The preceding lemmas show that if $\varphi\in{\mathcal L}^\alpha$, then,  for any $x\in\R^d$, $G_D(x,\cdot)\in \mathcal{H}$. Consequently,
the stochastic integral $\int_D G_D(x,y) d F(y)$ is well defined and
\begin{equation}\label{2.7}
E\bigg|\int_DG_D(x,y)\,dF(y)\bigg|^2\leq\|G_D(x,\cdot)\|_{L^\alpha(D)}^2 \le C_1^2.
\end{equation}
In addition, for each $p\in[1,\alpha']$, H\"older's inequality along with the hypercontractivity property and  (\ref{2.7}) yield
\begin{align}
&\bigg(E\bigg\|\int_D\,G_D(x,y)dF(y)\bigg\|_{L^{\alpha'}(D)}^{p}\bigg)^{1/p}
\leq
\bigg(E\bigg\|\int_D\,G_D(x,y)dF(y)\bigg\|_{L^{\alpha'}(D)}^{\alpha'}\bigg)^{1/\alpha'}\nonumber\\
&\quad
\leq
C\bigg(\int_D \bigg(E\bigg|\int_D\,G_D(x,y)dF(y)\bigg|^{2}\bigg)^{\alpha'/2}\,dx\bigg)^{1/\alpha'}\nonumber\\
&\quad\leq
C\bigg(\int_D \|G_D(x,\cdot)\|_{L^\alpha(D)}^{\alpha'}\,dx\bigg)^{1/\alpha'}
\leq C\,C_1. \label{2.8}
\end{align}
We next prove that the stochastic integral defines a H\"older continuous random field.

\begin{theorem}\label{t2.9} Fix $\lambda\in(0,1)$. Assume that $\varphi\in \mathcal{L}^\alpha(D)$ with
$\alpha$ in the interval $(1,\frac{d}{(d-2)(2-\lambda)\vee(d-1)\lambda})$. Then, 
for any $d\geq 4$, the Gaussian random field
$\{v(x)=\int_DG_D(x,y)\,dF(y), x\in D\}$ satisfies
\begin{equation}
\label{2.10}
E\left(|v(x)-v(z)|^2\right) \le C|x-z|^{2\lambda}.
\end{equation}
Therefore,  a.s. the sample paths are H\"older continuous  of order $\gamma\in(0,\lambda)$. 
\end{theorem}
\noindent{\sc Proof}. Fix $x,z\in D$. By the first inequality in (\ref{2.7}), 
$$
E\bigg|\int_D G_D(x,y)\,dF(y)-\int_D G_D(z,y)\,dF(y)\bigg|^2\le T_1+T_2,
$$
with
\begin{align*}
T_1&=\|G(x,\cdot)-G(z,\cdot)\|_{L^\alpha(D)}^2,\\
T_2&=\|(E_x(G(B_\tau,y))-E_z(G(B_\tau,y))\|_{L^\alpha(D)}^2.
\end{align*}
Schwarz inequality implies
\begin{align*}
T_1^{\frac{\alpha}{2}}&=\|G(x,\cdot)-G(z,\cdot)\|_{L^\alpha(D)}^\alpha
\le C\left(\int_D\bigg|\frac{1}{|x-y|^{d-2}}
-\frac{1}{|z-y|^{d-2}}\bigg|^{\alpha\lambda}\,dy\right)^{\frac{1}{2}}\\
&\quad\times \left(\int_D\bigg|\frac{1}{|x-y|^{d-2}}
+\frac{1}{|z-y|^{d-2}}\bigg|^{\alpha(2-\lambda)}\,dy\right)^{\frac{1}{2}}.
\end{align*}
The last integral is finite provided that $\alpha\in(0,\frac{d}{(d-2)(2-\lambda)})$. 

For the first one, we apply the mean value theorem to obtain
\begin{equation*}
\int_D\bigg|\frac{1}{|x-y|^{d-2}}
-\frac{1}{|z-y|^{d-2}}\bigg|^{\alpha\lambda}\,dy\le |x-z|^{\alpha\lambda}\int_D\bigg|\frac{|\xi|^{d-3}}{|x-y|^{d-2}|z-y|^{d-2}}\bigg|^{\alpha\lambda}\,dy,
\end{equation*}
with $\xi=\mu (x-y)+(1-\mu)(z-y)$ for some $\mu\in(0,1)$. 
Thus, if we define $D_1=D\cap\{|x-y|\leq |z-y|\}$ and $D_2=D\backslash D_1$, we have
\begin{eqnarray*}
\lefteqn{\int_D\bigg|\frac{|\xi|^{d-3}}{|x-y|^{d-2}|z-y|^{d-2}}\bigg|^{\alpha\lambda}\,dy}\\
&\leq&
\int_{D_1}\bigg|\frac{2|z-y|^{d-3}}{|x-y|^{d-2}|z-y|^{d-2}}\bigg|^{\alpha\lambda}\,dy
+\int_{D_2}\bigg|\frac{2|x-y|^{d-3}}{|x-y|^{d-2}|z-y|^{d-2}}\bigg|^{\alpha\lambda}\,dy\\
&\leq&C\Big(\int_{D_1}\frac{dy}{|x-y|^{\alpha\lambda(d-1)}}
+\int_{D_2}\frac{dy}{|z-y|^{\alpha\lambda(d-1)}}\Big)\\
&\leq&C\int_{D}\frac{dy}{|x-y|^{\alpha\lambda(d-1)}}\leq
C\int_0^{2}\frac{r^{d-1}\,dr}{r^{\alpha\lambda(d-1)}}.
\end{eqnarray*}
This last integral is finite if and only if $\alpha\in(0,\frac{d}{(d-1)\lambda})$. 
Summarising,  for $\alpha\in(1,\frac{d}{(d-2)(2-\lambda)\vee(d-1)\lambda})$,
\begin{equation}
\label{2.10'}
\Vert G(x,\cdot)-G(z,\cdot)\Vert_{L^\alpha(D)} \le C|x-z|^\lambda.
\end{equation}
For the second term, let us observe that by the strong Markov property,
$E_z(G(B_\tau,y))=E_x(G(B_\tau-x+z,y))$,
and therefore
\begin{eqnarray*}
T_2&=&C\|E_x(G(B_\tau,y))-E_x(G(B_\tau-x+z,y))\|_{L^\alpha(D)}^2\\
&\leq&C\bigg(\int_\Omega\|(G(u,\cdot))-G(u-x+z,\cdot)\|_{L^\alpha(D)}\,dP^x(u)\bigg)^2
\leq C|x-z|^{2\lambda},
\end{eqnarray*}
where in the last inequality we have applied (\ref{2.10'}). Thus, 
\begin{equation}\label{2.11}
\|G_D(x,\cdot)-G_D(z,\cdot)\|_{L^\alpha(D)}\leq C|x-z|^\lambda.
\end{equation}
Since the process $(v(x), x\in D)$ is Gaussian, the statement follows from  Kolmogorov's continuity criterion.

\cqd

\noindent{\bf Remark:} Assume $\varphi\in\mathcal{L}^\alpha$ for some $\alpha\in(1,\frac{d}{d-1})$. A slight
variant of the proof of Theorem \ref{t2.9} yields Lipschitz continuity of the sample paths of $(v(x), x\in D)$, a.s.

\subsection{Existence and uniqueness of the solution}\label{s2.2}

We give in this section a theorem on existence and uniqueness of solution for equation (\ref{1.1}) and we also establish some of the properties needed later.

We shall often use the next property proved in \cite{bu-pa:90}, Lemma 2.4; it is a consequence of the solvability of the
Dirichlet problem on $D$ and Poincar\'e's inequality (see \cite{gt}).

\noindent{\bf (P)} \quad There exists a constant $a>0$ such that for any $\varphi\in L^2(D)$,
\begin{equation}
\label{1.5}
\int_D \left(\int_D G_D(x,y)\varphi(y)dy\right) \varphi(x)dx \le -a\int_D\left(\int_D G_D(x,y)\varphi(y)dy\right)^2 dx.
\end{equation}
Unless otherwise stated, along this section $\lambda$ is a fixed number in $(0,1)$, $\alpha\in [1,\frac{d}{(d-2)(2-\lambda)\vee(d-1)\lambda})$, 
$\alpha'$ is its conjugate, that is, $\frac{1}{\alpha}+\frac{1}{\alpha'}=1$ and $a$ the positive constant given in (\ref{1.5}).

\begin{theorem}\label{t2.12}
Suppose that $g\in L^{\alpha'}(D)$, $f$ is of the form (\ref{1.2}) and satisfies (f1) and (f2) with a Lipschitz constant
 $L<\min(a, C_1^{-1})$.
Assume also that the correlation density $\varphi$ belongs to ${\mathcal L}^{\alpha}$. Then, there exists a unique
stochastic process solution to (\ref{1.1}).
\end{theorem}

\noindent{\sc Proof}. Let 
$${\mathcal B}=\big\{w:\, w\in{\mathcal C}(D),\ w|_{\partial D}=0\big\},$$
and consider the operator $T: {\mathcal B}\rightarrow {\mathcal B}$, defined by
$$T(w)(x)=w(x)-\int_DG_D(x,y)f(w(y))\,dy.$$
Clearly, ${\mathcal B}\subset L^{\alpha'}(D)$.

By H\"older inequality and (\ref{2.11}),
\begin{align*}
\left\vert \int_D\left(G_D(x,y)-G_D(z,y)\right) g(y)\right\vert&
\le \Vert g\Vert_{L^{\alpha'}(D)}\Vert G_D(x,y)-G_D(z,y)\Vert_{L^\alpha(D)}\\
&\le C|x-z|^\lambda.
\end{align*}
Together with Theorem \ref{t2.9}, this implies
$$b(x)=\int_DG_D(x,y)g(y)\,dy+\int_DG_D(x,y)\,dF(y)\in{\mathcal B},$$
for each $\omega$, a.s.

We next show that 
the operator equation
$Tw=b$ has a unique solution for any $b\in{\mathcal B}$, which happens if
$T$ is a bijective operator on ${\mathcal B}$. Uniqueness guarantees the
measurability of the process $(w(x), x\in D)$.

Let us first check that $T$ in one to one.
Fix $u$ and $v$ such that $Tu=Tv$. Then,
$u(x)-v(x)=\int_D G_D(x,y)(f(u(y))-f(v(y)))\,dy$. Multiplying both sides of
this equation by $f(u(x))-f(v(x))$  integrating over $D$ and using (\ref{1.5}) we obtain
\begin{eqnarray*}
\lefteqn{\int_D (u(x)-v(x))(f(u(x))-f(v(x)))\,dx}\\
&=&\int_D(f(u(x))-f(v(x)))\left(\int_D G_D(x,y)(f(u(y))-f(v(y)))\,dy\right)\,dx\\
&\leq&-a\int_D\bigg(\int_D G_D(x,y)(f(u(y))-f(v(y)))\,dy\bigg)^2\,dx\\
&=&-a\int_D(u(x)-v(x))^2\,dx.
\end{eqnarray*}
By (\ref{2.3})
$$\int_D (u(x)-v(x))(f(u(x))-f(v(x)))\,dx\geq -L\int_D(u(x)-v(x))^2\,dx.$$
Hence, 
$(a-L)\int_D(u(x)-v(x))^2\,dx\leq 0$. Since $L<a$ and $u,v\in L^{\alpha'}(D)$
with $\alpha'>2$, this implies $u(x)=v(x)$ for almost every $x\in D$.
\medskip

We next prove that $T$ is onto,  proceeding in a similar way than in \cite{catherine}. In the next arguments,
$\omega$ is a fixed element on a set of probability one.
\medskip

\noindent{\sf Step 1. A solution for a regular problem}. Let $b\in{\mathcal B}$, and 
$b_n\in{\mathcal C}^\infty_c(D)$, $n\ge 1$, such that $b_n\to b$ in $L^{\alpha'}(D)$. Clearly, 
the convergence also holds in $L^2(D)$. We construct  in the next Lemma \ref{l2.13} a
sequence of functions solving $Tu_n=b_n$ such that $u_n\to u$ in $L^2(D)$; the limit $u$ will be
our candidate for solution. 

Let us recall a basic result on the solution of nonlinear monotone operator equations (see e.g. 
\cite[Theorem 2.1, pg. 171]{lions},
or alternatively \cite[Theorem 26.A, pg. 557]{ze:90}):

\noindent{\bf (E)} Let $X$ be a reflexive Banach space; denote by $X^*$ its topological dual. Let $B: X\to X^*$ be 
a strictly monotone, coercive, hemicontinuous operator. Then, for any $k\in X^*$, the equation $Bw=k$ has a unique solution on $X$.
\smallskip

\begin{lemma}\label{l2.13}
For every $n\geq 1$, the boundary value problem
$$\Delta u_n(x)-f(u_n(x))=\Delta b_n(x) \mbox{ for } x\in D,\quad
u_n(x)=0 \mbox{ for }x\in \partial D,$$
has a unique solution $u_n(x)\in W_0^{1,2}(D)$.
\end{lemma}

\noindent{\sc Proof of Lemma \ref{l2.13}}. Set
$X=W_0^{1,2}(D)$, and define $B:W_0^{1,2}(D)\longrightarrow (W_0^{1,2}(D))^*$ by $Bu=-\Delta u+f(u)$;
that is, for any $w\in W_0^{1,2}(D)$,
$$
\langle Bu,w\rangle= \int_D\nabla w(x)\cdot\nabla u(x)\,dx+\int_Dw(y)f(u(y))\,dy
$$
The assumptions on $f$ imply that this operator satisfies the properties required in  {\bf (E)}. Moreover, 
 $k:=-\Delta b_n\in {\mathcal C}^\infty_c(D)\subset (W_0^{1,2}(D))^*$, for any $n\ge 1$. Thus,
the lemma follows from  {\bf (E)}.
\cqd
\bigskip

The sequence $\{u_n, n\ge 1\}$ satisfies
\begin{equation}\label{2.14}
u_n(x) = \int_D\,G_D(x,y)f(u_n(y))\,dy+b_n, \mbox{ for } x\in D, \quad u_n|_{\partial D}=0.
\end{equation}
Let us check that it is a Cauchy sequence in $L^2(D)$. By multiplying both
sides of the  equation satisfied by $u_n(x)-u_m(x)$ by $f(u_n(x))-f(u_m(x))$, integrating
over $D$, and owing to \eqref{2.3} and \eqref{1.5},
we obtain
\begin{eqnarray*}
-L\|u_n-u_m\|_{L^2(D)}^2&+&a\bigg\|\int_DG_D(x,y)(f(u_n(y))-f(u_m(y)))\,dy\bigg\|_{L^2(D)}^2\\
&\le&\int_D(b_n(x)-b_m(x))(f(u_n(x))-f(u_m(x)))\,dx.
\end{eqnarray*}
By using \eqref{2.14}, and the fact that each $u_n\in L^2(D)$, we obtain
\begin{eqnarray*}
&&(a-L)\|u_n-u_m\|_{L^2(D)}^2+a\|b_n-b_m\|_{L^2(D)}^2\\
&&\leq\int_D(b_n(x)-b_m(x))[(f(u_n(x))-f(u_m(x)))+2a(u_n(x)-u_m(x))]\,dx\\
&&\leq\|b_n-b_m\|_{L^2(D)}[2M+L\|u_n-u_m\|_{L^2(D)}+2a\|u_n-u_m\|_{L^2(D)}],
\end{eqnarray*}
that is,
$$ \|u_n-u_m\|_{L^2(D)}^2
\leq C\|b_n-b_m\|_{L^2(D)}[1+\|u_n-u_m\|_{L^2(D)}],$$
which implies that $\|u_n-u_m\|_{L^2(D)}\longrightarrow 0$ as
$n,m\to\infty$. 

Set $u=\lim_n u_n$ in $L^2(D)$.
\medskip

\noindent{\sf Step 2. $u$ is the solution}. We must prove
that $u\in \mathcal{B}$ and 
verifies \eqref{1.3}. That is, we would like to take limits in (\ref{2.14}).

We choose subsequences $u_n$ and $b_n$
(still denoted with the same subscripts) converging to $u$ and $b$
almost everywhere. We proceed in three steps.
\smallskip

\noindent{\sf Step 2.1.} Assume that $f$ is bounded (and continuous). Then
$$
u(x)-\int_DG_D(x,y)f(u(y))\,dy=b(x) \mbox{ for } x\in D,\quad u|_{\partial D}=0,
$$
by bounded convergence, and $u\in \mathcal{B}$.
\smallskip

\noindent{\sf Step 2.2.}  Assume that $f$ is bounded from below, that is, $f(x)\geq -N$ for every $x$
and some $N>0$. Set $f_n(x)=f_1(x)+(f_2(x)\wedge n)$, $n\ge 0$. 
Notice that each $f_n$ satisfies $(f1)$ and $(f2)$. Let
$$
u_n(x)=\int_DG_D(x,y)f_n(u_n(y))\,dy+b(x) \mbox{ for } x\in D, \quad u_n|_{\partial D}=0,
$$
be the solution constructed in {\sf Step 2.1}. 

We will need the following comparison statement. Its proof is very similar to that of Lemma 2.6 in \cite{bu-pa:90}
and therefore omitted.
\begin{lemma}\label{l2.15}
Let $f$ and $h$ satisfy $(f1)$, $(f2)$ and  $f(x)\geq h(x)$ for every
$x\in\R$. Let $b\in L^{\alpha'}(D)$ and
\begin{align*}
&u(x)-\int_DG_D(x,y)f(u(y))\,dy=b,\\
&v(x)-\int_DG_D(x,y)h(v(y))\,dy=b.
\end{align*}
Then, $u(x)\leq v(x)$ for almost every $x\in D$.
\end{lemma}
\smallskip

The sequence $(f_n, n\ge 0)$ is increasing; hence,
by Lemma \ref{l2.15}, the sequence of functions $(u_n, n\ge 0)$
satisfying
\begin{equation}
\label{2.16}
u_n(x)-\int_DG_D(x,y)f_1(u_n(y))\,dy-\int_DG_D(x,y)(f_2\wedge n)(u_n(y))\,dy=b(x),
\end{equation}
is decreasing. Set $u(x)=\inf_n u_n(x)$. Notice that
it is an a.s. finite function.
The function $f_1$ being bounded, we can take the limit inside the first integral in the left
hand-side of (\ref{2.16}). It remains to prove that we can also take the limit inside
the second integral. For this, we need some {\it a priori estimates} provided by the next two statements.

\begin{lemma}\label{l2.17}
The sequence $(u_n, n\ge 0)$ defined in (\ref{2.16}) satisfies
$$
\sup_n\|u_n\|_{L^{\alpha'}(D)}\leq \kappa,
$$
with $\kappa=\frac{(M+|f_2(0)|)C_1+\|b\|_{L^{\alpha'}(D)}}
{1-LC_1}$ and
$C_1= \sup_{x\in D}\Vert G_D(x,\cdot)\Vert_{L^{\alpha}(D)}$ (see Lemma \ref {l2.6}).
\end{lemma}

\noindent{\sc Proof}. Since $f$ is bounded from below,
$f_n$ is bounded by some constant depending on $n$ and 
thus each $u_n\in L^{\alpha'}(D)$.
Now, $|(f_2\wedge n)(u_n(x))|\leq |f_2(0)|+L|u_n(x)|$.
Thus, by H\"older inequality
$$
|u_n(x)|\leq (M+|f_2(0)|)C_1+LC_1\|u_n\|_{L^{\alpha'}(D)}+|b(x)|,
$$
and by integration on $D$, 
$$
\|u_n\|_{L^{\alpha'}(D)}\leq (M+|f_2(0)|)C_1+LC_1\|u_n\|_{L^{\alpha'}(D)}+\|b\|_{L^{\alpha'}(D)}.
$$
This yields the lemma.\cqd

\begin{corollary}\label{c2.18}
$\|u\|_{L^{\alpha'}(D)}\leq \kappa$, with $\kappa$ as in Lemma \ref{l2.17}.
\end{corollary}

\noindent{\sc Proof}. It is a consequence of Fatou's lemma. In fact,
\begin{equation*}
\|u\|_{L^{\alpha'}(D)}\le \liminf_n \|u_n\|_{L^{\alpha'}(D)}
\le \sup_n \|u_n\|_{L^{\alpha'}(D)}\le \kappa.
\end{equation*}

\cqd

Since $u_n(x)\downarrow u(x)>-\infty$ almost everywhere, we have
$|u(x)|\ind_{\{u>0\}}\leq |u_0(x)|\ind_{\{u>0\}}$,
and $|u(x)|\ind_{\{u<0\}}\leq \sup_n |u_n(x)|\ind_{\{u<0\}}$.
Thus,
\begin{eqnarray*}
|u_n(x)|
&=&|u_n(x)|\ind_{\{u>0\}}+|u_n(x)|\ind_{\{u<0\}}\\
&\leq&
|u_0(x)|\ind_{\{u>0\}}+|u(x)|\ind_{\{u<0\}}=\varphi(x),
\end{eqnarray*}
with $\varphi\in L^{\alpha'}(D)$ and therefore,
$$
|(f_2\wedge n)(u_n(x))|\leq|f_2(u_n(x))|\leq M+|f_2(0)|+L\varphi(x)\in L^{\alpha'}(D).
$$
Thus, by the dominated convergence theorem and the continuity of $f_2$,
$$
\lim_{n\to\infty}\int_DG_D(x,y)(f_2\wedge n)(u_n(y))\,dy=
\int_DG_D(x,y)f_2(u(y))\,dy.
$$
Summarizing, if $f$ is bounded from below, there exists $u$ satisfying
$$
u(x)=\int_DG_D(x,y)f(u(y))\,dy+b(x) \mbox{ for } x\in D,  \quad u |_{\partial D}=0,
$$
and $\|u\|_{L^{\alpha'}(D)}\leq \kappa$.
\smallskip

\noindent{\sf Step 2.3. $f$ satisfies (f1) and (f2)}. Set $f_n=f_1+(f_2\vee (-n))$, $n\ge 0$. By the results obtained in the previous step,
there exists $u_n$ such that
$$
u_n(x)=\int_DG_D(x,y)f_n(u_n(y))\,dy+b(x)\mbox{ for } x\in D,\qquad \quad u |_{\partial D}=0,
$$
and
$\sup_n\|u_n\|_{L^{\alpha'}(D)}\leq \kappa$.

The sequence $(f_n, n\ge 0)$ is decreasing; hence, by Lemma \ref{l2.15}
$\{u_n, n\ge 0\}$ is increasing. Set $u(x)=\sup_n u_n(x)$ for a.e. $x$.
As in {\sf Step 2.2},  it suffices to prove that
$$
\lim_{n\to\infty}\int_DG_D(x,y)(f_2\vee(-n))(u_n(y))\,dy=
\int_DG_D(x,y)f_2(u(y))\,dy.
$$
Lemma \ref{l2.17} and Fatou's Lemma implies $\Vert u\Vert_{L^{\alpha'}(D)}<\kappa$.

Since $u_n(x)\uparrow u(x)$, we have $|u(x)|\ind_{\{u>0\}}= \sup_n|u_n(x)|\ind_{\{u>0\}}$,
and $|u(x)|\ind_{\{u<0\}}\leq |u_0(x)|\ind_{\{u<0\}}$. Thus, 
\begin{eqnarray*}
|u_n(x)|
&=&|u_n(x)|\ind_{\{u>0\}}+|u_n(x)|\ind_{\{u<0\}}\\
&\leq&
|u(x)|\ind_{\{u>0\}}+|u_0(x)|\ind_{\{u<0\}}=\psi(x),
\end{eqnarray*}
with $\psi\in L^{\alpha'}(D)$.\\
Observe that  $|f_2\vee (-n)|\leq |f_2|$. Therefore
$$
|(f_2\vee( -n))(u_n(y))|\leq|f_2(u_n(y))|\leq M+|f_2(0)|+L\psi(x).
$$
Thus, by the bounded convergence theorem and the continuity of $f_2$,
$$
\lim_{n\to\infty}\int_DG_D(x,y)(f_2\vee (-n))(u_n(y))\,dy=
\int_DG_D(x,y)f_2(u(y))\,dy.
$$
Hence, we have proved the existence of $u\in L^{\alpha'}(D)$ satisfying
\begin{equation}
\label{2.19}
u(x)=\int_DG_D(x,y)f(u(y))\,dy+b(x) \mbox{ for } x\in D,\end{equation}
$u|\partial D=0$. The terms in the right hand-side of (\ref{2.19}) belong to  $\mathcal{B}$; therefore, so does $u$.

\cqd

\subsection{Properties of the solution}\label{s2.3}

The solution of (\ref{1.1}) given in Theorem \ref{t2.12} possesses some important properties, as follows.

\begin{lemma}\label{l2.20}
With the same assumptions as in Theorem \ref{t2.9},  for any $p\in[1,\alpha']$, 
$$
\|u\|_{L^p(\Omega; L^{\alpha'}(D))}\leq C.
$$
\end{lemma}

\noindent{\sc Proof}.
The solution $u$ to equation \eqref{1.1}
satisfies $\|u\|_{L^{\alpha'}(D)}\leq \kappa$, where $\kappa$ is as in Lemma \ref{l2.17}. 
Hence, it suffices to check that $b\in L^p(\Omega;L^{\alpha'}(D))$, for any  $p\in[1,\alpha']$, where
$$b(x)=\int_DG_D(x,y)g(y)\,dy+\int_DG_D(x,y)dF(y).$$
The function $x\mapsto \int_DG_D(x,y)g(y)\,dy$ is continuous and deterministic; consequently it belongs to
$L^p(\Omega;L^{\alpha'}(D))$. By virtue of (\ref{2.8}) the same conclusion holds for the stochastic
integral $\int_DG_D(\cdot,y)dF(y)$.
\cqd

\begin{theorem}\label{t2.21} With the same hypotheses as in Theorem \ref{t2.9}, 
for any $p\in[1,\alpha']$, the solution 
 $u$ to equation \eqref{1.1}  satisfies 
 \begin{equation}
 \label{2}
 E\left(|u(x)-u(z)|\right)^p\leq C|x-z|^{p\lambda}.
 \end{equation}
 Consequently, a.s. the sample paths are $\gamma$-H\"older continuous with $\gamma\in(0,\lambda-\frac{d}{\alpha'})$.
\end{theorem}

\noindent{\sc Proof}.
Clearly, 
$\big(E(|u(x)-u(z)|^p\big)^{1/p}\leq \sum_{i=1}^3I(x,z)$,
with
\begin{align*}
I_1(x,z)&=
\bigg(E\bigg|\int_D(G_D(x,y)-G_D(z,y))f(u(y))\,dy\bigg|^p\bigg)^{1/p}\\
I_2(x,z)&=\bigg(E\bigg|\int_D(G_D(x,y)-G_D(z,y))g(y)\,dy\bigg|^p\bigg)^{1/p}\\
I_3(x,z)&=\bigg(E\bigg|\int_D(G_D(x,y)-G_D(z,y))\,dF(y)\bigg|^p\bigg)^{1/p}.
\end{align*}
By H\"older inequality, Lemma \ref{l2.20} and (\ref{2.11}),
we have
\begin{align*}
I_1(x,z)&\leq\|G_D(x,\cdot)-G_D(z,\cdot)\|_{L^\alpha(D)}\big(M+|f_2(0)|+LE(\|u\|_{L^{\alpha'}(D)}^p)^{1/p}\big)\\
&\leq C|x-z|^\lambda.
\end{align*}
A similar bound holds for the term $I_2(x,z)$. 

The hypercontractivity property and (\ref{2.10}) yield $I_3(x,z)\le C |x-z|^\lambda$.
Thus, we have proved (\ref{2}). The statement about H\"older continuity follows from Kolmogorov's criterion.
\cqd

\section{Numerical approximations in $L^2(D)$--norm}\label{s3}

This section is devoted to give a finite differences scheme for the spde (\ref{1.1}). We have shown that in dimension $d\ge 4$, the Green function
does not belong to $L^2(D)$. Thus, the method used in \cite{GM} for $d\le 3$ does not apply directly. Instead, we combine a smoothing of $G_D$ with a weighted discretization.

\subsection{Construction of the smoothing $G_D^\varepsilon$}\label{s3.1}

We start by introducing some technical background (see \cite{folland2}). Let ${\mathbb T}^d$  be the 
$d$-dimensional torus endowed with the Haar measure $dx$.
Any function
$f:[-1,1)^d\longrightarrow {\R}$ is identified with  $F:{\mathbb T}^d\longrightarrow {\R}$
defined by 
$$F(e^{i\pi x}):=F(e^{i\pi x_1},\dots,e^{i\pi x_d}),$$
that is, $F(e^{i\pi x})=f(x)$. Then
$$\int_{(-1,1)^d}f(t)\,dt=\int_{{\mathbb T}^d}F(e^{i\pi t})\,dt$$
Consider the odd extension of  $G_D(x,y)$  in the $y$-variables, that is, if $y_j\in(0,1)$, $j=1,\dots,d$, we define
$$
G_D(y_1,\dots,-y_i,\dots,y_d)=-G_D(y_1,\dots,y_i,\dots,y_d).
$$
The new function, still denoted by $G_D(x,\cdot)$, is now defined on $(-1,1)^d$. 
Let
$$
{\mathbb G}_D^{x}(e^{i\pi y})=G_D(x,y)
$$
its identification in the torus setting. Observe that ${\mathbb G}_D^{x}$ satisfies
$$
\int_{{\mathbb T}^d}|{\mathbb G}_D^{x}(e^{i\pi y})|^\alpha\,dy
=\int_{(-1,1)^d}|G_D(x,y)|^\alpha\,dy=2^d\|G_d(x,\cdot)\|_{L^\alpha(D)}^\alpha.
$$

Let $\psi(x)\in{\mathcal C}^\infty_c((-1,1))$ be an even function,  $0\leq \psi\leq 1$,
$\int_{-1}^1\psi=1$ and $\Psi(x)=\prod_{i=1}^d\psi(x_i)$.
Clearly, $\Psi(x)\in{\mathcal C}^\infty_c((-1,1)^d)$ and it is an even function in each variable $x_i$.
Set
$$
\Phi(e^{i\pi x}) := \prod_{i=1}^d\phi(e^{i\pi x_i}) := \prod_{i=1}^d\psi(x_i)=\Psi(x).
$$
The functions
$$\Phi_\varepsilon(e^{i\pi x}):=\frac1{\varepsilon^d}\Psi\big(\frac x{\varepsilon}\big):=\Psi_\varepsilon(x),$$
$\varepsilon>0$, define an approximation of the identity in ${\mathbb T}^d$.

By means of $\Phi_\varepsilon(e^{i\pi \cdot})$, we define the smoothing of $G_D(x,\cdot)$ as follows:
$$
{\mathbb G}_D^{x,\varepsilon}(e^{i\pi y})=
\int_{{\mathbb T}^d}{\mathbb G}_D^{x}(e^{i\pi (y-u)})\Phi_\varepsilon(e^{i\pi u})\,du.
$$
Clearly, ${\mathbb G}_D^{x}(e^{i\pi y})\in L^\alpha({\mathbb T})$, therefore
${\mathbb G}_D^{x,\varepsilon}\to {\mathbb G}_D^{x}$ in $L^\alpha(\tor)$ as $\varepsilon\to 0$.
Define
$$G_D^\varepsilon(x,y)={\mathbb G}_D^{x,\varepsilon}(e^{i\pi y}).$$
It is not difficult to check that, for any $\tilde y=(y_1,\dots,-y_i,\dots,y_d)$, $i=1,\dots,d$, we have
$G_D^\varepsilon(x,\tilde y)=-G_D^\varepsilon(x, y)$, that is, $G_D^\varepsilon(x, y)$ is odd in all the variables $y_i$.
In addition, for any $\alpha\in[1,\frac{d}{d-2})$,
\begin{equation}\label{3.1}
\sup_\varepsilon\sup_{x\in D}\|G_D^\varepsilon(x,\cdot)\|_{L^\alpha(D)}\leq CC_1,
\end{equation}
with the same constant $C_1$ defined in Lemma \ref{l2.6}.

The next result provides a bound for the error of the smoothing.

\begin{lemma}\label{l3.2} Fix $\lambda\in(0,1)$ and as in Theorem \ref{t2.9}, assume that $\varphi\in \mathcal{L}^\alpha(D)$ with
$\alpha\in(0,\frac{d}{(d-2)(2-\lambda)\vee(d-1)\lambda})$. 
There exists a constant $C$ such that, for every $\varepsilon>0$,
$$\|G_D-G_D^\varepsilon\|_{L^\alpha(D\times D)}\leq C\varepsilon^\lambda.$$
\end{lemma}

\noindent{\sc Proof}. Since $G_D(x,\cdot)$, $G_D^\varepsilon(x,\cdot)$
are odd in the $y$-variables, we have
\begin{eqnarray*}
\|G_D(x,\cdot)-G_D^\varepsilon(x,\cdot)\|_{L^\alpha(D)}
&=&
2^{-d}\|G_D(x,\cdot)-G_D^\varepsilon(x,\cdot)\|_{L^\alpha((-1,1)^d)}\\
&=&
2^{-d}\|\G_D^x-\G_D^{x,\varepsilon}\|_{L^\alpha(\tor^d)}.
\end{eqnarray*}
Thus, by Fubini's theorem
\begin{equation*}
\|G_D-G_D^\varepsilon\|_{L^\alpha(D\times D)}^\alpha=
2^{-d\alpha}\big\|\,\|\G_D^x(\cdot)-\G_D^{x,\varepsilon}(\cdot)\|_{L^\alpha(D)}^\alpha\,\big\|_{L^\alpha(\tor^d)}^\alpha,
\end{equation*}
with the variable $x$ integrated over $D$. 
H\"older's inequality with respect to the finite measure on $\tor^d$ given by $\Phi_\varepsilon(u)du$
implies, for any fixed $y$,
\begin{align*}
&\|\G_D^x(e^{i\pi y})-\G_D^{x,\varepsilon}(e^{i\pi y})\|_{L^\alpha(D)}^\alpha\\
&\quad=
\bigg\|\G_D^x(e^{i\pi y})-\int_{\tor^d} \G_D^x(e^{i\pi (y-u)})\Phi_\varepsilon(u)\,du\bigg\|_{L^\alpha(D)}^\alpha\\
&\quad=
\bigg\|\int_{\tor^d}(\G_D^x(e^{i\pi y})- \G_D^x(e^{i\pi (y-u)}))\Phi_\varepsilon(u)\,du\bigg\|_{L^\alpha(D)}^\alpha\\
&\quad \leq
\int_{\tor^d}\big\|\G_D^x(e^{i\pi y})- \G_D^x(e^{i\pi (y-u)})\big\|_{L^\alpha(D)}^\alpha\Phi_\varepsilon(u)\,du.
\end{align*}
By \eqref{2.11} and Fubini's theorem, we see that
$$ \|G_D-G_D^\varepsilon \|_{L^\alpha(D\times D)}^\alpha \le C \int_{\tor^d}|u|^{\alpha \lambda}\Phi_\varepsilon(u)\,du.$$
Taking into account that $\Phi_\varepsilon (u)\neq 0$ if and only if $u\in(-\varepsilon,\varepsilon)^d$,
we conclude the proof.
\cqd

Next, we give the Fourier expansion of $G_D^\varepsilon(x,y)$.
Let us recall that for any $d\geq 1$, the set of functions 
$$v_{\beta}(x)=\sin(\beta_1\pi x_1)\cdots\sin(\beta_d\pi x_d),\ \beta\in I^d,$$
is an orthogonal complete system in $L^2(D)$, with $\|v_{\beta}\|_{L^2(D)}=2^{-d/2}$.  
  
\begin{lemma}\label{l3.4}
For any $\varepsilon>0$, we have
$$
G_D^\varepsilon(x,y)=\sum_{\beta\in I^d}\frac{-\hat\Psi(\varepsilon\alpha)2^{d}}{\pi^2|\alpha|^2}
v_{\beta}(x)v_{\beta}(y),
$$
in $L^2(D\times D)$ and a.e., where $\hat\Psi(\xi)$ denotes the Fourier transform of $\Psi$.
\end{lemma}

In the sequel, we shall often use the following remark: The function  $\hat\Psi$ 
is rapidly decreasing, that is,
for any multiindex $\theta\in (0,\infty)^d$, there exists a constant
$C(\theta)$ such that $|\xi|^\theta |\hat\Psi(\xi)|\leq C(\theta)$.
\smallskip

\noindent{\sc Proof of Lemma \ref{l3.4}}.
By virtue of Young's inequality for convolutions
\begin{align}
\|G_D^{\varepsilon}(x,\cdot)\|_{L^2(D)}
&=
2^{-d}\|\G_D^{x,\varepsilon}(e^{i\pi \cdot})\|_{L^2(\tor^d)}\nonumber\\
&\leq
2^{-d}\|\G_D^x\|_{L^\alpha(\tor^d)}\|\Phi_\varepsilon\|_{L^\rho(\tor^d)}\\
&=\|G_D(x,\cdot)\|_{L^\alpha(D)}\|\Psi_\varepsilon\|_{L^\rho((-1,1)^d)},\label{bound0}
\end{align}
with $\frac{1}{2}+1=\frac{1}{\alpha}+\frac{1}{\rho}$.
For any $\rho\ge 1$, $\sup_\varepsilon \Vert\Psi_\varepsilon\Vert_{L^\rho((-1,1)^d)}<C$. Hence, by Lemma \ref{l2.6},
\begin{equation}
\label{bound1}
\sup_\varepsilon\sup_{x\in D}\|G_D^{\varepsilon}(x,\cdot)\|_{L^2(D)}\le C_1,
\end{equation}
and consequently $G_D^\varepsilon\in L^2(D\times D)$.

We next compute the Fourier coefficients of  $G_D^\varepsilon$, as follows. Set $v_{\beta}(y):=\\
\sigma_\beta(e^{i\pi y})$. 
By Fubini's theorem
\begin{align*}
\int_DG_D^\varepsilon(x,y)v_{\beta}(y)\,dy
&=2^{-d}\int_{(-1,1)^d}G_D^\varepsilon(x,y)v_{\beta}(y)\,dy\\
&=2^{-d}\int_{\tor^d}\G_D^{x,\varepsilon}(e^{i\pi y})\sigma_\beta(e^{i\pi y})\,dy\\
&=2^{-d}\int_{\tor^d}\left(\int_{\tor^d}\G_D^{x}(e^{i\pi (y-u)})
\Phi_\varepsilon(e^{i\pi u})\,du \right) \sigma_\beta(e^{i\pi y})\,dy\\
&=2^{-d}\int_{\tor^d}\left(\int_{\tor^d}\G_D^{x}(e^{i\pi (y-u)})\sigma_\beta(e^{i\pi y})\,dy\right)
\Phi_\varepsilon(e^{i\pi u})\,du,\\
\end{align*}
Set
$K(e^{i\pi u})= \int_{\tor^d}\G_D^{x}(e^{i\pi (y-u)})\sigma_\beta(e^{i\pi y})\,dy$.
By a change of variables,
$$K(e^{i\pi u})
=\int_{\tor^d}\G_D^{x}(e^{i\pi z})\sigma_\beta(e^{i\pi (z+u)})\,dz
=\int_{(-1,1)^d}G_D(x,z)v_{\beta}(z+u)\,dz.
$$
We compute this last integral using the formula
\begin{align*}
v_{\beta}(z+u)&=\prod_{i=1}^d\sin \beta_i\pi(z_i+u_i)\\
&=
\prod_{i=1}^d(\sin\beta_i\pi z_i\cos\beta_i\pi u_i+\cos\beta_i\pi z_i\sin\beta_i\pi u_i).
\end{align*}
The product on the right-hand-side of the former inequality consists of
a sum of terms, each of them being a product of $d$ factors either of the form
$\sin\beta_i\pi z_i\cos\beta_i\pi u_i$ or $\cos\beta_i\pi z_i\sin\beta_i\pi u_i$.
Since $G_D(x,z)$ is odd in all the $z$-variables, integrals of terms containing
factors $\cos\alpha_i\pi z_i\sin\alpha_i\pi u_i$
will be zero. Thus, 
\begin{align*}
K(e^{i\pi u})
&=\int_{(-1,1)^d}G_D(x,z)\prod_{i=1}^d(\sin\beta_i\pi z_i\cos\beta_i\pi u_i)\,dz\\
&=\bigg(\prod_{i=1}^d\cos\beta_i\pi u_i\bigg)
\int_{(-1,1)^d}G_D(x,z)v_{\beta}(z)\,dz\\
&=\bigg(\prod_{i=1}^d\cos\beta_i\pi u_i\bigg)
2^d\int_{D}G_D(x,z)v_{\beta}(z)\,dz\\
&=\bigg(\prod_{i=1}^d\cos\beta_i\pi u_i\bigg)\frac{-2^d}{\pi^2|\beta|^2}v_{\beta}(x),
\end{align*}
where in the last equality we have used the properties of $G_D$.
Therefore, 
\begin{align*}
&\int_DG_D^\varepsilon(x,y)v_{\beta}(y)\,dy
=2^{-d}\int_{\tor^d}K(e^{i\pi u})\Phi_\varepsilon(e^{i\pi u})\,du,\\
&\qquad \quad=2^{-d}\int_{(-1,1)^d}
\bigg(\prod_{i=1}^d\cos\beta_i\pi u_i\bigg)\frac{-2^d}{\pi^2|\beta|^2}v_{\beta}(x)
\Psi_\varepsilon(u)\,du\\
&\qquad \quad=\frac{-1}{\pi^2|\beta|^2}v_{\beta}(x)
\prod_{i=1}^d\int_{(-1,1)}\cos\beta_i\pi x\psi_\varepsilon(x)\,dx,
\end{align*}
But
\begin{equation*}
\int_{(-1,1)}\cos\beta_i\pi x\psi_\varepsilon(x)\,dx=
\int_{(-\varepsilon,\varepsilon)}e^{-i\beta_i\pi x}
\frac1\varepsilon\psi\big(\frac x\varepsilon\big)\, dx=
\hat \psi(\varepsilon\alpha_i),
\end{equation*}
where $\hat \psi$ stands for the  Fourier transform of $\psi$ in $\R$. Therefore,
\begin{equation*}
\int_DG_D^\varepsilon(x,y)v_{\beta}(y)\,dy
=\frac{-1}{\pi^2|\beta|^2}v_{\beta}(x)
\prod_{i=1}^d\hat \psi(\varepsilon\beta_i)
=\frac{-1}{\pi^2|\beta|^2}v_{\beta}(x)
\hat \Psi(\varepsilon\beta).
\end{equation*}
This finishes the proof of the lemma.
\cqd

\subsection{Construction of the discrete approximations}\label{s3.2} Let us introduce some notation.
For any $n\ge 1$, we consider the grid of $\bar D=[0,1]^d$ given by 
$$
\mathcal{G}=\Big\{\frac{j}{n}=\bigg(\frac{j_1}{n},\dots, \frac{j_d}{n}\bigg):\ j_k=0,1\dots,n,\ k=1,
\dots,d\Big\}\subset \bar D.
$$
For any point $\frac{j}{n}\in\mathcal{G}$, we set 
$
D_j=\bigg[\frac{j_1}{n},\frac{j_1+1}{n}\bigg)\times\dots\times
\bigg[\frac{j_d}{n}, \frac{j_d+1}{n}\bigg)
$
and define $\kappa_n(x)=\frac{j}{n}$ for each $x\in D_j$. 

On the space  $X=\big\{ u:\ u=\{u_i\}_{i\in
I^d_n}\big\}=\R^{(n-1)^d}$ endowed with the Hilbert-Schmidt
norm, we consider the second order difference operator $A:X\to X$ defined  by
$$
(Au)_i=\sum_{j=1}^d n^2\big[u_{i-e_j}-2u_i+u_{i+e_j}\big],$$
where $\{e_j\}_{j=1}^d$ is the canonical basis of $\R^d$,
The set of vectors of $X$,
$$ \left(
\bigg(\frac{2}{n}\bigg)^{d/2}U_\beta, \beta\in I^d_n\right),\, 
(U_\beta)_i=v_{\beta}\big(\frac{i}{n}\big), \,  i\in I^d_n $$
is an orthonormal system in $X$ of eigenvectors of $A$,
with eigenvalues
$$\lambda_\beta=-\pi^2(\beta_1^2c_{\beta_1}+\dots+\beta_d^2c_{\beta_d}),$$ where
$c_l=\sin^2\bigg(\frac{l\pi}{2n}\bigg)\left(\frac{l\pi}{2n}\right)^{-2}$.
Notice that $\frac{4}{\pi^2}\leq c_l\leq 1$ for every $1\leq l\leq n-1$.

In the sequel, we consider  the lexicographic order in $I^d_n$. Denoting by $U$ the $(n-1)^d$ matrix
whose rows are the vectors $U_{\beta_j}$, (here $\beta_j$, $j=1,\cdots,(n-1)^d$, denotes the lexicographic enumeration of $I_n^d$) we have 
$$A=U^tDU,$$
with $D$ the square diagonal matrix with entries $D_{j,j}=\lambda_{\beta_j}$.


For any $\varepsilon>0$, define $D^{\epsilon}$ the square diagonal matrix in dimension $(n-1)^d$ with
diagonal elements
$$\lambda^\epsilon_{\beta_j}=\frac{\lambda_{\beta_j}}{\hat\Psi(\varepsilon\beta_j)}.$$

We also consider a sequence $(g_n, n\ge 1)$ of step functions defined on $D$ such that 
for $n$ big enough, $\|g-g_n\|_{L^{\alpha'}(D)}^2\leq C/n^2$ for some positive
constant $C$.

The discrete approximations of $u$ are defined first on points of $\mathcal{G}$ as follows.
If $\frac{j}{n}\in \mathcal{G}\cap\partial D$, define $u_n^\varepsilon(\frac{j}{n})=0$ 
(boundary conditions). For $\frac{j}{n}$, with $j\in I_n^d$, we define $u_n^\varepsilon(\frac{j}{n})$ to be the solution
of the system
\begin{equation}
\label{3.6}
(U^t D^\varepsilon U) u_n^\varepsilon=f(u_n^\varepsilon)+g_n+n^d{\mathbf F},
\end{equation}
where ${\mathbf F}$ is
the vector $(F(D_i), i\in I_n^d)$, with $F(D_i)= \int_{\R^d}\ind_{D_i} dF(y)$, and
$g_n(x)=g_n(\kappa_n(x))$, $n\ge 1$.
Finally, for any $x\in D$ we define $u_n^\varepsilon(x)=u_n^\varepsilon(\kappa_n(x))$. 


We prove in Proposition \ref{p3.9} that a solution to equation (\ref{3.6}) exists. Moreover, proceeding as in \cite{GM}, 
it is easy to check that $u_n^\varepsilon$ satisfies the {\it mild} equation
\begin{align}
u_n^\varepsilon(x)&=\int_DG_{D,n}^\varepsilon(x,y)f(u_n^\varepsilon(y))\,dy
+\int_DG_{D,n}^\varepsilon(x,y)g_n(y)\,dy\\
&+\int_DG_{D,n}^\varepsilon(x,y)dF(y),\label{3.7}
\end{align}
with
\begin{equation}\label{3.8}
G_{D,n}^\varepsilon(x,y)=\sum_{\beta\in I^d_n}\frac{\hat\Psi(\varepsilon\beta)2^d}
{\lambda_{\beta}}v_{\beta}(\kappa_n(x))v_{\beta}(\kappa_n(y)).
\end{equation}
We will prove later that an appropriate sequence $u_n:=u_n^{\varepsilon(n)}$
of these  approximations converges to the solution
of \eqref{1.1} in the space $L^p(\Omega;L^2(D))$, for any $p\in[1,\alpha']$, with an specific rate.
%

\begin{proposition}\label{p3.9} Assume that $\varphi\in \mathcal{L}^\alpha$ for some $\alpha\in[1,\frac{d}{d-2})$. Suppose 
also that $g\in L^{\alpha'}(D)$, f is of the form (\ref{1.2}) and satisfies $(f1)$ and $(f2)$ with a Lipschitz constant $L<4d$. Then
the system (\ref{3.6}) has a unique solution.
\end{proposition}
\noindent{\sc Proof}. We apply the classical result quoted in {\bf (E)} of Section \ref{s2.2} to $X=X^*=\R^{(n-1)^d}$, endowed with the Hilbert-Schmidt
norm, $Bu=fu-A^\varepsilon u$, $A^\varepsilon=U^t D^\varepsilon U$, $f(u)=\{f(u_i)\}$ and $b\in X$ with components
$b_i=-g_n(\frac{i}{n})-n^d{\mathbf F}(D_i)$. The set of
vectors of $X$,
$ \Big\{\bigg(\frac{2}{n}\bigg)^{d/2}U_{\beta_i}, i=1, \cdots, (n-1)^d\Big\}$,
$(U_{\beta_i})_j=\varphi_{\beta_i}\bigg(\frac{{\mathbf j}}{n}\bigg)$,
where $i=1,\dots,(n-1)^d$ and ${\mathbf j}$ is the $j$-th vector of
$I^d_n$ ordered with the lexicographic order, is an orthonormal system
of $X$ of eigenvectors of $A^\varepsilon$ with eigenvalues
\begin{equation}\label{3.10}
\frac{\lambda_{\beta}}{\hat\Psi(\varepsilon\beta)}=-\frac{\pi^2(\beta_1^2c_{\beta_1}+\dots+\beta_d^2c_{\beta_d})}
{\hat\Psi(\varepsilon\beta)},\quad
c_l=\frac{\sin^2\bigg(\frac{l\pi}{2n}\bigg)}{\bigg(\frac{l\pi}{2n}\bigg)^2},
\end{equation}
satisfying $\frac{4}{\pi^2}\leq c_l\leq 1$ for every $1\leq l\leq n-1$. 
Thus, $4|\beta|^2\leq|\lambda_\beta|\leq \pi^2|\beta|^2$.
This property together with the obvious bound
 $|\hat\Psi(\varepsilon\beta)|\leq\|\Psi\|_{L^1(\R^d)}=1$, imply
$-\lambda_{\beta}/\hat\Psi(\varepsilon\beta)\geq 4|\beta|^2\geq 4d$.

The operator $B$ is strictly monotone. Indeed, for any $u,v\in X$,
$u\neq v$, we write $u-v=\sum_{\beta\in I^d_n}(u_\beta-v_\beta)\big(\frac{2}{n}\big)^{d/2}U_\beta$. Then,
\begin{eqnarray*}
\lefteqn{\langle Bu-Bv, u-v\rangle}\\
&=&
\langle-\sum_{\beta\in I^d_n}\frac{\lambda_\beta}{\hat\Psi(\varepsilon\beta)}(u_\beta-v_\beta)
\bigg(\frac{2}{n}\bigg)^{d/2}U_\alpha+(f(u)-f(v)), u-v\rangle\\
&=&
\sum_{\beta\in I^d_n}\frac{\pi^2(\beta_1^2c_{\beta_1}+\dots+
\beta_d^2c_{\beta_d})}{\hat\Psi(\varepsilon\beta)}
(u_\beta-v_\beta)^2\\
& &+\sum_{i\in I^d_n}(u_i-v_i)(f(u_i)-f(v_i))\\
&\geq& 4d|u-v|^2-L|u-v|^2=(4d-L)|u-v|^2,
\end{eqnarray*}
where we have applied (\ref{2.3}).

Let us now prove that $B$ is coercive. Since $B$ is strictly monotone, and $B(0)=f(0)$, we have
\begin{eqnarray*}
\langle Bu,u\rangle
&=&
\langle B(u)-B(0), u\rangle+\langle B(0),u\rangle\geq (4d-L)|u|^2
-|\langle B(0),u\rangle|\\
&\geq& (4d-L)|u|^2-|f(0)|\,|u|.
\end{eqnarray*}
Therefore,
$$
\frac{\langle Bu,u\rangle}{|u|}\geq(4d-L)|u|-|f(0)|,
$$
which implies $\lim_{|u|\to\infty}\frac{\langle Bu,u\rangle}{|u|}=\infty$.

Finally, since for any $u,v\in X$, the functions $t\mapsto U^t D^\varepsilon U(u+tv)$ and $t\mapsto f(u+tv)$ are continuous, 
$B$ is an hemicontinuous operator.
 
 \cqd

\subsection{Properties of the regularized and truncated kernels}\label{s3.3}

This section is devoted to prove some integrability properties of the kernels  $G_{D,n}^\varepsilon$ and estimates of
the discrepancy between $G_D^\varepsilon$ and $G_D$, $G_{D,n}^\varepsilon$, respectively.

Along the section, $\theta$ is a fixed positive real number satisfying $\theta>2d-4$.

\begin{lemma}\label{l3.11}
Let $\varepsilon:=\varepsilon(n)=n^{(2d-4-\theta)/\theta}$. 
There exists a constant $C(\theta)>0$ such that
\begin{equation}
\label{bound2}
\sup_{n\ge 1}\sup_{x\in D}\|G_{D,n}^{\varepsilon(n)}(x,\cdot)\|_{L^2(D)}\leq C(\theta).
\end{equation}
\end{lemma}
\noindent{\sc Proof}. 
The system
$\{v_{\beta}(\kappa_n(y))\}$ is orthogonal in $L^2(D)$. Moreover,\\
 $\sup_{x\in D}|v_{\beta}(x)|\le C$ and 
$|\lambda_\beta|\geq 4|\beta|^2$. Thus,
$$
\|G_{D,n}^\varepsilon(x,\cdot)\|_{L^2(D)}
=
C\sum_{\beta\in I^d_n}\frac{\hat\Psi(\varepsilon\beta)^2}{\lambda_{\beta}^2}
\leq
\frac{C(\theta)}{\varepsilon^\theta}
\sum_{\beta\in I^d_n}\frac{1}{|\beta|^{4+\theta}}.
$$
for any $\theta>0$.

The sum in the right-hand-side of this expression is comparable with a  Riemann sum for the integral
of the function $|x|^{-4-\theta}$ on a region away from the origin. Observe that $\sqrt d\leq |\beta|
\leq n\sqrt d$, for $\beta\in I^d_n$, and $|x|^{-4-\theta}$ is radial and decreasing. Denoting by
$Q_\beta=\prod_{i=1}^d(\beta_i-1,\beta_i)$, we have $|\beta|\geq |x|$, and
\begin{eqnarray*}
\frac1{n^d}\sum_{\beta\in I^d_n}\frac{1}{|\beta|^{4+\theta}}
&\leq&
2^d\int_{(1/2,1)^d}|x|^{-4-\theta}\,dx+\sum_{\beta\in I^d_n}\int_{Q_\beta}
|x|^{-4-\theta}\,dx\\
&\leq&
C\int_{1/2}^{n\sqrt d}\frac{r^{d-1}\,dr}{r^{4+\theta}}=Cr^{d-4-\theta}\big|_{1/2}^{n\sqrt d}
\leq
C n^{d-4-\theta}.
\end{eqnarray*}
Thus,
\begin{eqnarray*}
\|G_{D,n}^\varepsilon(x,\cdot)\|_{L^2(D)}
\leq
\frac{C(\theta)}{\varepsilon^\theta}n^dn^{d-4-\theta}=
\frac{C(\theta)}{\varepsilon^\theta}n^{2d-4-\theta}.
\end{eqnarray*}
Hence, choosing $\varepsilon=n^{(2d-4-\theta)/\theta}$
we obtain (\ref{bound2}).

\cqd


The next result gives an estimate in the $L^2(D\times D)$--norm of the approximation of the smoothed Green function by truncation. First, we recall
the following facts (see Lemma 3.2 \cite{GM}): For any $\beta\in I^d$, $x,z\in \R^d$,
\begin{align}
&|v_{\beta}(x)-v_{\beta}(z)|\leq C|\beta|\,|x-z|,\label{bound6}\\
&\Big|\frac{-1}{\pi^2|\beta|^2}-\frac{1}{\lambda_\beta}\Big|\leq\frac{C}{|\beta|n}.\label{bound7}
\end{align}

\begin{lemma}\label{l3.13}
Set $\varepsilon(n)=n^{(2d-4-\theta)/\theta}$. For
every $\gamma>0$ there exists a constant $C(\gamma)$ such that
\begin{align}
\label{bound10}
\sup_{x\in D}\|G_D^{\varepsilon(n)}(x,\cdot)-G_{D,n}^{\varepsilon(n)}(x,\cdot)\|_{L^2(D)}\leq C(\gamma)n^{-\gamma}.
\end{align}
\end{lemma}

\noindent{\sc Proof}. 
For simplicity, we write $\varepsilon$ instead of $\varepsilon(n)$. 
By the definitions of the kernels $G_D^\varepsilon$ and $G_{D,n}^\varepsilon$, we can write
\begin{equation*}
\int_D |G_D^\varepsilon(x,y)-G_{D,n}^\varepsilon(x,y)|^2\,dy\le 2 (A(x)+B(x)),
\end{equation*}
with
\begin{align*}
A(x)&=\int_D\Big|\sum_{\beta\in I^d\backslash I^d_n}\frac{-2^d\hat\Psi(\varepsilon\beta)}{\pi^2|\beta|^2} v_{\beta}(x)v_{\beta}(y)\Big|^2\,dy\\
B(x)&=\int_D\bigg| \sum_{\beta\in I^d_n}
\bigg[\frac{-2^d\hat\Psi(\varepsilon\beta)}{\pi^2|\beta|^2}v_{\beta}(x)v_{\beta}(y)\\
&\quad -
\frac{2^d\hat\Psi(\varepsilon\beta)}{\lambda_\beta}v_{\beta}(\kappa_n(x))v_{\beta}(\kappa_n(y))\bigg]\bigg|^2\,dy
\end{align*}
Since $\hat\Psi(\varepsilon\beta)$ is a rapidly decreasing function, owing to the
orthogonality of the system defined by the functions $v_{\beta}$,
we have
\begin{align*}
A(x)&\le\frac{C(b)}{\varepsilon^b}\sum_{\beta\in I^d\backslash I^d_n}\frac{1}{|\beta|^{4+b}}
\leq
\frac{C(b)n^d}{\varepsilon^b}\int_{n\sqrt d}^\infty\frac{r^{d-1}\,dr}{r^{4+b}}\\
&\leq
\frac{C(b)}{\varepsilon^b}n^{2d-4-b}=C(b) n^{(2d-4)(1-\frac{b}{\theta})},
\end{align*}
for every $b>d-4$, where $C(b)$ is a positive constant depending on $b$.

Fix $\gamma >0$ and  choose $b>2d-4>d-4$ such that $(2d-4)\big(1-\frac b\theta\big)=-2\gamma$.
We obtain $\sup_{x\in D}A(x)\le C(\gamma) n^{-2\gamma}$.

Clearly, $B(x)\le C(B_1(x)+B_2(x))$, with
\begin{align*}
B_1(x)&=\int_D\bigg|\sum_{\beta\in I^d_n}\frac{\hat\Psi(\varepsilon\beta)}{|\beta|^2}
(v_{\beta}(x)v_{\beta}(y)-v_{\beta}(\kappa_n(x))v_{\beta}(\kappa_n(y)))\bigg|^2\,dy,\\
B_2(x)&=\int_D\bigg|\sum_{\beta\in I^d_n}\bigg(-\frac{1}{\lambda_\beta}-\frac1{\pi^2|\beta|^2}\bigg)
\hat\Psi(\varepsilon\beta)
v_{\beta}(\kappa_n(x))v_{\beta}(\kappa_n(y)))\bigg|^2\,dy.
\end{align*}
By virtue of (\ref{bound6}) 
\begin{align*}
B_1^{\frac{1}{2}}(x)\le &C\sum_{\beta\in I^d_n}\frac{|\hat\Psi(\varepsilon\beta)|}{|\beta|^2}\frac{|\beta|}{n} \le \frac{C(l)n^d}{n\varepsilon^l}\int_{1/2}^{n\sqrt d}\frac{r^{d-1}\,dr}{r^{1+l}}\\
&\le\frac{C(l)n^{d-1}}{\varepsilon^l}n^{d-l-1}=\frac{C(l)n^{2d-l-2}}{\varepsilon^l},
\end{align*}
for any $l>0$.

Let $l$ be such that  $2d-2-(2d-4)\frac l\theta=-\gamma$;
we obtain $\sup_{x\in D}B_1^{\frac{1}{2}}(x)\le C(\gamma) n^{-\gamma}$.

The orthogonality of the vectors $v_{\beta}(\kappa_n(x))$, (\ref{bound7}) and the properties of $\Psi$ imply
\begin{align*}
B_2(x)&\le C \sum_{\beta\in
I^d_n}\bigg|\frac{1}{\lambda_\beta}+\frac1{\pi^2|\beta|^2}\bigg|^2|\hat\Psi(\varepsilon\beta)|^2
\leq
C\sum_{\beta\in I^d_n}\frac{C(k)}{\varepsilon^k|\beta|^k}\frac{1}{|\beta|^2n^2}\\
&\le \frac{C(k)n^d}{n^2\varepsilon^k}\int_{1/2}^{n\sqrt d}\frac{r^{d-1}\,dr}{r^{2+k}}
\leq
\frac{C(k)n^d}{n^2\varepsilon^k}n^{d-k-2},
\end{align*}
for any $k>0$. Choosing $k$ such that $(2d-4)\big(1-\frac k\theta\big)=-2\gamma$, we obtain
$\sup_{x\in D}B_2(x)\le C(\gamma) n^{-2\gamma}$.

We have thus finished the proof.

\cqd


\begin{lemma}
\label{l3.16}
With the same assumptions as in Lemma \ref{l3.2}, set $\varepsilon(n)=\\
n^{(2d-4-\theta)/\theta}$.
There exists a constant $C:=C(\lambda,\theta,d)$ not depending on $n$, such that
\begin{equation}
\label{bound12}
\|G_D-G_{D,n}^\varepsilon\|_{L^\alpha(D\times D)}\leq Cn^{-\gamma},
\end{equation}
with $\gamma=\frac{\lambda}{\theta}(\theta+4-2d)$.
\end{lemma}

\noindent{\sc Proof}. It is a consequence of Lemmas \ref{l3.2} and \ref{l3.13}.
\cqd
\smallskip

\noindent{\bf Remark}  In Lemma \ref{l3.16}, the parameter $\theta$ can be chosen arbitrarily large. Therefore, for any given $\delta >0$, one could obtain
(\ref{bound12}) with $\gamma = \lambda-\delta$. On the other hand, we do not have an explicit control on the dependence of the constant $C$ on
$\theta$. Actually, this constant appears in the formulation of the rapidly decreasing property of the Fourier transform of the regularising kernel $\Psi$.


\subsection{Properties of the approximations}\label{s3.4}

One of the consequences of the properties of the modified Green kernels established in the preceding section
is the following {\it a priori} estimate for the solution of equation (\ref{3.7}). This is the main result.

\begin{proposition}\label{p3.18} Fix $\theta>2d-4$. Assume that the Lipschitz constant in $(f2)$ satisfies $L\le \min(4d, [C(\theta)]^{-1})$, with $C(\theta)$ given in (\ref{bound2}), and the
 hypotheses of Proposition \ref{p3.9}. Let 
 $\varepsilon(n)=n^{(2d-4-\theta)/\theta)}$. Then, for any  $p\in[1,\alpha']$$$
\sup_{n\geq 1}
\left(\|u^{\varepsilon(n)}_n\|_{L^p(\Omega;L^{\alpha'}(D))}\right)\leq \tilde C(\theta),
$$
for some positive constant $\tilde C(\theta)$ depending on $\theta$.
\end{proposition}

\noindent{\sc Proof}. Let us write $\varepsilon$ instead of 
$\varepsilon(n)$. Since $\alpha\le 2$, H\"older's inequality and (\ref{bound2}) yield
$\sup_n\sup_{x\in D} \Vert G_{D,n}^\varepsilon\Vert_{L^{\alpha}(D)}\le C$. Then, H\"older's inequality
and the properties on $f$ imply
\begin{align*}
|u^\varepsilon_n(x)|&\le \left|\int_DG_{D,n}^\varepsilon(x,y)dF(y)\right| + \sup_n\sup_{x\in D} \Vert G_{D,n}^\varepsilon\Vert_{L^{\alpha}(D)}\\
&\quad\times\left(\sup_n\Vert g_n\Vert_{L^{\alpha'}(D)}+M+f(0)+L\|u^\varepsilon_n\|_{L^{\alpha'}(D)}\right).
\end{align*}
Notice that since $u_n^\varepsilon$ is a step function, its $L^{\alpha'}$--norm is finite. 

Proceeding as in (\ref{2.8}) with the Green function $G_D$ replaced by $G_{D,n}^\varepsilon$ and using
(\ref{bound2}), we see that for some positive constant $\hat{C}(\theta)$,
\begin{equation*}
E\left(\left\Vert \int_D G_{D,n}^\varepsilon(\cdot,y) dF(y)\right\Vert_{L^{\alpha'}(D)}^p\right)\le \hat C(\theta).
\end{equation*}
Since $L\le \left[C(\theta)\right]^{-1}$, the announced result follows easily.

\cqd

\subsection{Convergence results}\label{s3.6} 

We devote this section to the proof of the approximation of the solution of (\ref{1.1}) by means of the discretized scheme defined in (\ref{3.7}),
with an appropriate choice of the smoothing parameter $\varepsilon$. Here is the statement.

\begin{theorem}\label{t3.19} Fix $\lambda\in(0,1)$ and $\theta>2d-4$. We assume that the hypotheses of Theorem \ref{t2.12}
and Propositions \ref{p3.9} and \ref{p3.18} are satisfied. That is, $\varphi\in \mathcal{L}^\alpha$ with $\alpha\in\left(0,\frac{d}{(d-2)(2-\lambda)\vee(d-1)\lambda}\right)$ and $L\le \min(a,4d,C_1^{-1}, [C(\theta)]^{-1})$.
Set $\gamma=\frac{\lambda}{\theta}(\theta+4-2d)$. Then, there exist a constant $C(\gamma)>0$ such that for any $p\in [1,\alpha']$,
\begin{equation}
\label{aprox1}
\|u-u^{\varepsilon(n)}_n\|_{L^p(\Omega; L^2(D))}\leq C(\gamma) n^{-\gamma\frac{\alpha}{2\alpha'}},
\end{equation}
where $\varepsilon(n)=n^{(2d-4-\theta)/\theta)}$, $n\ge 1$. 
\end{theorem}

\noindent{\sc Proof}. Throughout this proof we write $\varepsilon$ instead of $\varepsilon(n)$, for the sake of simplicity.
Set 
\begin{align}
T(x)&=\int_D[G_D(x,y)-G_{D,n}^\varepsilon(x,y)]f(u^\varepsilon_n(y))\,dy\nonumber\\
& +\int_DG_D(x,y)[g(y)-g_n(y)]\,dy
+\int_D[G_D(x,y)-G_{D,n}^\varepsilon(x,y)]g_n(y)\,dy\nonumber\\
& +\int_D[G_D(x,y)-G_{D,n}^\varepsilon(x,y)]\,dF(y),\label{3.16}
\end{align}
so that
\begin{equation}\label{3.17}
u(x)-u^\varepsilon_n(x)=\int_DG_D(x,y)[f(u(y))-f(u^\varepsilon_n(y))]\,dy
+T(x).
\end{equation}
We multiply both sides of (\ref{3.17}) by $f(u(x))-f(u^\varepsilon_n(x))$, then we apply the inequality (\ref{2.3}) and integrate over $D$ and 
apply (\ref{1.5}).
We obtain 
\begin{align*}
\lefteqn{-L\|u-u^\varepsilon_n\|_{L^2(D)}^2}\\
&\leq
\int_D\Big(\int_DG_D(x,y)[f(u(y))-f(u^\varepsilon_n(y))]\,dy\Big)\big(f(u(x))-f(u^\varepsilon_n(x))\big)\,dx\\
& +\int_DT(x)\big(f(u(x))-f(u^\varepsilon_n(x))\big)\,dx\\
&\le -a\Big\|\int_DG_D(x,y)[f(u(y))-f(u^\varepsilon_n(y))]\,dy\Big\|_{L^2(D)}^2\\
& +\int_DT(x)\big(f(u(x))-f(u^\varepsilon_n(x))\big)\,dx.
\end{align*}

Now, by (\ref{3.17}),
\begin{eqnarray*}
\lefteqn{\Big\|\int_DG_D(x,y)[f(u(y))-f(u^\varepsilon_n(y))]\,dy\Big\|_{L^2(D)}^2}\\
&=&
\|u-u^\varepsilon_n\|_{L^2(D)}^2+\|T\|_{L^{\alpha}(D)}^2-2\int_D(u(x)-u^\varepsilon_n(x))T(x)\,dx.
\end{eqnarray*}
Thus,
\begin{eqnarray*}
\lefteqn{-L\|u-u^\varepsilon_n\|_{L^2(D)}^2
\leq
-a\|u-u^\varepsilon_n\|_{L^2(D)}^2-a\|T\|_{L^{\alpha}(D)}^2}\\
& &+2a\int_D(u(x)-u^\varepsilon_n(x))T(x)\,dx
+\int_D\big(f(u(x))-f(u^\varepsilon_n(x))\big) T(x)\,dx.
\end{eqnarray*}
From here, applying H\"older inequality and the hypothesis $(f1)$ and $(f2)$
\begin{align}
\lefteqn{(a-L)\|u-u^\varepsilon_n\|_{L^2(D)}^2}\nonumber\\
&\leq 2a\int_D(u(x)-u^\varepsilon_n(x))T(x)\,dx
+\int_D \big(f(u(x))-f(u^\varepsilon_n(x))\big)T(x)\,dx\nonumber\\
&\leq
2a\|u-u^\varepsilon_n\|_{L^{\alpha'}(D)}\|T\|_{L^\alpha(D)}
+\|f(u)-f(u_n^\varepsilon)\|_{L^{\alpha'}(D)}\|T\|_{L^\alpha(D)}\nonumber\\
&\leq
2a\|u-u^\varepsilon_n\|_{L^{\alpha'}(D)}\|T\|_{L^{\alpha}(D)}
+\|T\|_{L^{\alpha}(D)}(2M+L\|u-u^\varepsilon_n\|_{L^{\alpha'}(D)})\nonumber\\
&=
(2a+L)\|u-u^\varepsilon_n\|_{L^{\alpha'}(D)}\|T\|_{L^{\alpha}(D)}
+2M\|T\|_{L^{\alpha}(D)}. \label{3.18}
\end{align}

Notice that, by virtue of Lemma \ref{l2.20} and Proposition \ref{p3.18}, for the choice $\varepsilon:=\varepsilon(n)$ and for any $1\le q\le \alpha'$,
\begin{equation}
\label{3.19}
\sup_nE(||u-u_n^\varepsilon||^q_{L^{\alpha'}(D)}) \le C.
\end{equation}
Assume $p\in[1, 2]$.  Applying H\"older's inequality, then (\ref{3.18}),  Schwarz inequality and (\ref{3.19}) with $q=2$ yields 
\begin{equation}
\label{3.20}
\bigg(E(\|u-u^\varepsilon_n\|_{L^2(D)}^p\bigg)^{1/p}
\leq C\bigg(E\big(\|T\|_{L^{\alpha}(D)}^2\big)\bigg)^{1/4}.
\end{equation}

Assume now $p\in]2,\alpha']$. Owing to (\ref{3.18}), Schwarz inequality and (\ref{3.19})
\begin{align}
\label{3.21}
\Big(E\big(\|u-u^\varepsilon_n\|_{L^2(D)}^p\big)\Big)^{2/p}
\leq
C\Big(E\bigg(\|T\|_{L^{\alpha}(D)}^p\bigg)\Big)^{1/p},
\end{align}
Hence, we have reduced the proof of the theorem to that of giving estimates of the $L^p(\Omega)$--norm of $\Vert T\Vert_{L^{\alpha}(D)}$, for $p\in[2,\alpha']$.
 
By the definition of $T(x)$ (see (\ref{3.16})) and 
H\"older's inequality, for any $\alpha>1$ we have 
\begin{align*}
\Vert T\Vert_{L^{\alpha}(D)}&\le \Vert G_D-G_{D,n}^\varepsilon\Vert_{L^\alpha(D\times D)} \Vert f(u_n^\varepsilon)\Vert_{L^{\alpha'}(D)}\\
& + \Vert G_D\Vert_{L^\alpha(D\times D)}\Vert g-g_n\Vert_{L^{\alpha'}(D)}+\Vert G_D-G_{D,n}^\varepsilon\Vert_{L^\alpha(D\times D)} \Vert g_n\Vert_{L^{\alpha'}(D)}\\
&+\Big(\int_D\bigg|\int_D[G_D(x,y)-G_{D,n}^\varepsilon(x,y)]\,dF(y)\bigg|^\alpha\,dx\Big)^{1/\alpha}.
\end{align*}
Lemmas \ref{l3.16}, \ref{l2.6} and  the hypotheses on the coefficients $f$ and $g$ yield 
\begin{align}
\|T\|_{L^{\alpha}(D)}
&\leq C(\gamma) n^{-\gamma}\nonumber\\
&\quad\times \bigg(M+1+
L\bigg(\int_D|u^\varepsilon_n(y)|^{\alpha'}\,dy\bigg)^{1/\alpha'}+|f_2(0)|+\|g\|_{L^{\alpha'}(D)}\bigg)\nonumber\\
&+C_1 n^{-1}+
 \Vert \int_D (G_D(\cdot,y)-G_{D,n}^\varepsilon(\cdot,y) dF(y)\Vert_{L^{\alpha'}(D)}, \label{3.22}
\end{align}
where $\gamma=\frac{\lambda}{\theta}(\theta+4-2d)$.

Let $p\in[1,\alpha']$.  According to (\ref{2.8}), where we replace $G_D$ by $G_D-G_{D,n}^\varepsilon$, we obtain
\begin{align*}
&\bigg(E\big(\Vert \int_D (G_D(\cdot,y)-G_{D,n}^\varepsilon(\cdot,y) dF(y)\Vert_{L^{\alpha'}(D)}^p\big)\bigg)^{1/p}\\
&\quad \le C\bigg(\int_D\|G_D(x,\cdot)-G_{D,n}^\varepsilon(x,\cdot)\|_{L^\alpha(D)}^{\alpha'}\,
dx\bigg)^{1/\alpha'}\\
&=C \bigg(\int_D\|G_D(x,\cdot)-G_{D,n}^\varepsilon(x,\cdot)\|_{L^\alpha(D)}^{\alpha}
\|G_D(x,\cdot)-G_{D,n}^\varepsilon(x,\cdot)\|_{L^\alpha(D)}^{\alpha'-\alpha}\,dx\bigg)^{1/\alpha'}\\
&\le C\|G_D(x,\cdot)-G_{D,n}^\varepsilon(x,\cdot)\|_{L^\alpha(D\times D)}^{\alpha/\alpha'},
\end{align*}
by Lemma \ref{l2.6} and (\ref{bound2}).
Together with  (\ref{3.22}), Proposition \ref{p3.18}, and Lemma \ref{l3.16}, this implies 
\begin{equation*}
\bigg(E\big(\|T\|_{L^{\alpha}(D)}^p\big)\bigg)^{1/p}
\leq C(\gamma)\left( n^{-\gamma}+n^{-\gamma\frac{\alpha}{\alpha'}}\right) \le C(\gamma)n^{-\gamma\frac{\alpha}{\alpha'}}.
\end{equation*}
Substituting this bound in (\ref{3.20}), (\ref{3.21}) yields the upper estimate (\ref{aprox1}) and therefore, the theorem is proved.

\cqd
\smallskip

\noindent{\bf Remark:} In the previous theorem,  the values of $p$ cannot be arbitrarily large; therefore, one cannot obtain
a rate of convergence for the sample paths, as is the case for $d\le 3$ (see for instance Corollary 2.4 in \cite{GM}).
\smallskip

\noindent{\bf Acknowledgement:} Part of this work has been done during a visit of the first named author to the
Institut de Matem\`atica, Universitat de Barcelona.

\end{document}